\documentclass{amsart}

\usepackage{amscd}
\usepackage[mathscr]{eucal}
\usepackage{enumerate}
\newtheorem{theorem}{Theorem}[section]
\newtheorem{lemma}[theorem]{Lemma}
\newtheorem{proposition}[theorem]{Proposition}

\theoremstyle{remark}

\newtheorem{corollary}[theorem]{Corollary}
\theoremstyle{definition}
\newtheorem{definition}[theorem]{Definition}

\theoremstyle{remark}
\newtheorem{remark}[theorem]{Remark}

\theoremstyle{remark}

\theoremstyle{definition}

\theoremstyle{remark}

\newenvironment{myequation}{{\addtocounter{theorem}{1}}\begin{equation}}{\end{equation}}

\newcommand{\rr}{{\mathbf  R}}

\newcommand{\ra}{\rightarrow}

\newcommand{\mymbox}[1]{{\mbox{~ #1 ~}}}

\newcommand{\fs}{\,\,\mbox{for some}\,\,}

\newcommand{\st}{\,\,\mbox{such that}\,\,}

\newcommand{\lb}{{\mathcal{O}}}
\newcommand{\mult}{{\mathrm {mult}}}

\newcommand{\supp}{\mathrm{Supp}}
\begin{document}
\title[A characterization of double covers of curves]{A characterization of double covers of curves in terms of the ample cone of second symmetric product}
%\subtitle{Do you have a subtitle?\\ If so, write it here}
\author{Kungho Chan}
%
%\institute{{Kungho Chan} at KIAS Hoegiro 87(207-43 Cheongnyangni
%2-dong),Dongdaemun-gu,Seoul 130-722,Korea}
%
\date{Received: date / Revised version: date}
% The correct dates will be entered by the editor
%
\maketitle
\begin{abstract}
We investigate the nef cone spanned by the diagonal and the fibre
classes of second symmetric product of a curve of genus $g$. This
$2$-dimensional nef cone gives a characterization of double covers
of curves of genus $\le \frac{g-1}{8}$. This is a generalization of
a result by Debarre \cite[Proposition 8]{De}.

% give the Mathematics Subject Classification here if any, e.g.:
% \subclass{00A00}
\end{abstract}
% If there is no abstract please say so using
% \noabstract instead
%
\section{Introduction} \label{intro}

Let $C$ be a smooth integral curve of genus $g \ge 2$, $C^{(2)}$ the
second symmetric product of $C$, and $x$ and $\delta$ the fibre and
the diagonal numerical classes respectively in the
N{\'{e}}ron-Severi space ${\mathrm{N}}^1(C^{(2)})_{\rr}$ where the
fibre class is the numerical class of the curve $C+a$ for any fixed
point $a \in C$. We are interested in the restriction of the nef
cone ${\mathrm{Nef}}(C^{(2)})$ on the plane spanned by $x$ and
$\delta$. It has two boundaries, one is
\[(g-1)x+(\delta/2)\] which is orthogonal to the diagonal class in the
cone of effective curves on $C^{(2)}$, and the other is
\[(\tau(C)+1)x-(\delta/2)\] where \[\tau(C) := \inf \{\mu \ge 0~|~(\mu
+ 1)x-({\delta  / 2})~\mbox{is~nef~on~}C^{(2)}\}.\] Since the
self-intersection of $(\tau(C)+1)x - ({\delta  / 2})$ is
non-negative, \[\tau(C) \ge \sqrt{g}.\]

When $C$ is very general, ${\mathrm{N}}^1(C^{(2)})_{\rr}$ is
generated by $x$ and $\delta$, and it was conjectured that $\tau(C)
= \sqrt{g}$ for $g \ge 4$ and $g$ not a perfect square \cite[Section
3]{CK} \cite[Conjecture 1.1]{Ro}.

On the other side, it is interesting to know what the value of
$\tau(C)$ is when $C$ is special. The first answer was given by
Lazarsfeld \cite[Example 1.5.13]{L1}. Among other things, the same
result was also given by Debarre \cite{De} and independently by Kong
\cite{Kon}. For $g \ge 3$,

\[\tau(C) = g \mbox{~if and only if~} C \mbox{~is hyperellipic.}\]
Actually, Debarre proved more in the same paper.

\begin{proposition} \cite[Proposition 8]{De} \label{debarre:result2}
Let $C$ be a smooth projective curve of genus $g \ge 2$.

\begin{enumerate}[(a)]
\item If $C$ hyperelliptic, then $\tau(C) = g$.
\item If $C$ is non-hyperelliptic, then
\begin{itemize}
\item If $g = 3$, then $\tau(C) = \frac{9}{5}$.
\item If $g = 4$, then $\tau(C) = 2$.
\item If $g \ge 5$, then $\tau(C) \le g-2$ where the equality holds if and
only if $C$ is bielliptic.
\end{itemize}
\end{enumerate}
\end{proposition}

It characterizes the biellipticity of $C$ with the value of $\tau(C)
= g-2$.

Along this direction, we give the main theorem in this paper.

\begin{theorem} \label{main:result}
Fix an integer $k \ge 0$. If $C$ is a smooth integral curve of genus
$g
> \max\{2k+1,4k-3\}$, then

\begin{enumerate}
\item[]  \begin{equation*}
\begin{split}{}
 &\tau(C)  \ge g-k \mbox{~if~and~only~if~}\\
 &\mbox{there~exists~a~smooth~integral~curve~} H \\
 &\mbox{of~genus~} q
\mbox{~with~} q \le \frac{k}{2} \st C \mbox{~is~a~double~cover~of~}
H.\\
\end{split}
 \end{equation*}

\item[] Furthermore in this case $H$ is unique up to isomorphic and
\[\tau(C) = g-2q.\]
\end{enumerate}
\end{theorem}

For $g \ge 6$, choose the greatest integer $k$ such that
$\frac{g+3}{4}> k \ge 2$. Then $k \ge \frac{g-1}{4}$. The following
is an immediate consequence.

\begin{corollary} \label{immed_conseq}
Suppose $C$ is a smooth irreducible curve of genus $g \ge 6$. We
have
\begin{itemize}
\item[(a)]  If $\tau(C) \ge \frac{3g+1}{4}$, then

\[C {~is~a~double~cover~of~a~curve~of~genus~} \frac{g-\tau(C)}{2}.\]

\item[(b)] If $C$ is a double cover of curve of genus $q \le
\frac{g-1}{8}$, then

\[\tau(C) = g-2q.\]

\end{itemize}
\end{corollary}

More concisely, it gives a characterization of double covers of
curves as following.
\begin{corollary}
Let ${\mathcal{M}}_g$ be the set of smooth integral curves of genus
$g \ge 6$.\\Then,
\[\{ C \in {\mathcal{M}}_g | \tau(C) \ge \frac{3g+1}{4}\} = \{ C \in
{\mathcal{M}}_g | C \mbox{~is a double cover of curve of genus~}q
\le \frac{g-1}{8}\}.\]
\end{corollary}

In the last section, we compute the values of $\tau(C)$ for some
other classes of special curves.

{\it{Acknowledgement.}} {This paper is part of my PhD thesis. I
thank my advisor Lawrence Ein for getting me interested in this
question and helping me with comments and suggestions.}

\section{Curves induced by Pencils} \label{pencil:curve}
With the notations in the previous section, for computing the
constant $\tau(C)$, we need to create some curves on the second
symmetric products. We make use of the pencils for this purpose. In
this section, we derive some basic properties about how these pencil
curves being related to the value of $\tau(C)$.

On ${\mathrm{N}}^1(C^{(2)})_{\rr}$, we have the following
intersection equations
\begin{myequation} \label{products}
(x.x) = 1, (x.\delta) = 2, (\delta.\delta) = 4-4g.
\end{myequation}

For convenience of computations, for any divisor $D$ on $C^{(2)}$,
we can always find unique $n$ and $\gamma$ such that

\[D \equiv (n+\gamma)x - \gamma (\delta/2) + \xi\] where $\xi \in
<x,\delta>^{\bot}$. We set up the notation
\[ (A,B)+\xi := (A+B)x - B (\delta/2) + \xi.\]
For divisors
\[D_1 \equiv (n_1,\gamma_1) + \xi_1 \mbox{~and~} D_2 \equiv
(n_2,\gamma_2) + \xi_2,\] by \eqref{products}, we have
\begin{myequation} \label{product:eqt}
\begin{split}
(D_1.D_2) &= n_1 n_2 - \gamma_1 \gamma_2 g + (\xi_1 . \xi_2),\\
D_1 \pm D_2 &\equiv (n_1 \pm n_2, \gamma_1 \pm \gamma_2) + (\xi_1
\pm
\xi_2).\\
\end{split}
\end{myequation}

\begin{definition}
Suppose $C$ admits a base-point-free pencil $g^1_d$ of degree $d$.
Then one can define a reduced curve on $C^{(2)}$ as
\[\Gamma(g^1_d) := \{ a+b \in C^{(2)}~| g^1_d-a-b \ge 0 \}.\]
\end{definition}
The numerical class of $\Gamma(g^1_d)$ in
${\mathrm{N}}^1(C^{(2)})_{\rr}$ is
\begin{myequation}\label{pencil:class}
d\cdot x - (\delta /2).
\end{myequation}

\begin{lemma} \label{basic:bounds}
We have the following simple facts on the value of $\tau(C)$.
\begin{enumerate}[(a)]
\item $\tau(C) \ge \sqrt{g}$.
\item For any curve $C^{(2)} \supset \gamma \equiv (n,\gamma) + \xi$, $n$ must be a positive integer.

\item Assume $C$ admits a base-point-free pencil $g^1_d$ of degree
$d$. Then, \begin{itemize}
\item $\tau(C) \ge \frac{g}{d-1}$.
\item If $d < 1+\sqrt{g}$, then $\tau(C) >
\sqrt{g}$.
\item If $d \ge 1 + \sqrt{g}$ and the curve $\Gamma(g^1_d)$ is
irreducible, then $\tau(C) \le d-1$.
\end{itemize}
\item For any smooth curve $C$, $\tau(C) \le g$.
\end{enumerate}
\end{lemma}

\begin{proof}
\begin{enumerate}[(a)]
\item Since $(\tau(C),1)$ is nef, $(\tau(C),1).(\tau(C),1) = \tau(C)^2
-g \ge 0$.

\item For any curve $C^{(2)} \supset \gamma \equiv (n,\gamma) + \xi$, since $x$ is ample,
$(x.C) = n$ which must be a positive integer.
\item
\begin{itemize}
\item From \eqref{pencil:class} we know $\Gamma(g^1_d) \equiv (d-1,1)$ and
\[(\tau(C),1).(d-1,1) = \tau(C)(d-1) - g \ge 0\] by the definition of
$\tau(C)$.

\item Trivial consequence from the last result.
\item Since $\Gamma(g^1_d)^2 = (d-1,1)^2 = (d-1)^2 -g \ge 0$, then $\Gamma(g^1_d) \equiv
(d-1,1)$ is nef, and hence $\tau(C) \le d-1$.
\end{itemize}
\item It is well-known that any $C$ always admits a
base-point-free pencil $A$ of degree $g+1$ and the $\Gamma(A)$ is
irreducible (\cite{BEL} Lemma 1.3 and Lemma 1.4 (i)), then $\tau(C)
\le g$ from above.
\end{enumerate}
\end{proof}

In \cite{Ko}, Kouvidakis proved the following.
\begin{proposition}\label{kouvidakis}
Suppose $C$ is a smooth integral projective curve and $C$ admits a
base-point-free pencil $g^1_d$ of degree $d$. If $d \le 1+
[\sqrt{g}]$ and $\Gamma(g^1_d)$ is irreducible, then
\[\tau(C) = \frac{g}{d-1}.\]
\end{proposition}

Kouvidakis's result means that if a smooth integral curve $C$ has a
low degree base-point-free pencil and the curve induced by the
pencil on $C^{(2)}$ is irreducible, then the class of
$\Gamma(g^1_d)$ is orthogonal to the boundary
$(\tau(C)+1)x-(\delta/2)$ and it computes $\tau(C)$.

Actually, this fact is also true even when the induced curve is
reducible. To this end, we have to consider the situation of an
arbitrary effective divisor on $C^{(2)}$.

First, we need the following modified version of a lemma by Ross.
\begin{lemma} [\cite{Ro} Lemma 2.2] \label{in:support}
Suppose $C$ is a smooth integral curve of genus $g \ge 2$. We have,
\begin{enumerate}[(a)]
\item If $\tau(C) > \sqrt{g}$, then there exists an irreducible
curve $\Gamma$ on $C^{(2)}$ such that $\Gamma \equiv (n,\gamma) +
\xi$ with $\xi \in <x,\delta>^{\bot}$, $\gamma\sqrt{g} > n$ and
$\tau(C) = \frac{\gamma g}{n}$ (i.e. this curve $\Gamma$ computes
$\tau(C)$).
\item If there exists an effective divisor $D$ on $C^{(2)}$ such
that $D \equiv (a,b)$ with integers $a,b>0$, then either $\tau(C)
\le a/b,$ or an irreducible component of $\supp(D)$ computes
$\tau(C)$. In the latter situation,

\[ \tau(C) = \max\limits_{C' \subset \supp (D)} \{ \mu ~|~ \mu = R(C') \}\] where the maximum is taken over all
irreducible components of $\supp(D)$ and

\[R(C') := \frac{(\delta.C') - 2(x.C')}{2(x.C')}.\]

In particular, if $C$ admits a base-point-free pencil $g^1_e$ and
$\tau(C) > e-1$, then

\[ \tau(C) = \max\limits_{C' \subset \Gamma(g^1_e)} \{ R(C') \}\] where the maximum is taken over all
irreducible components of $\Gamma(g^1_e)$.

\end{enumerate}
\end{lemma}
\begin{proof}
For (a). Since $(\tau(C) + 1)x - (\delta/2)$ is nef and not ample
and since $(\tau(C),1)^2 > 0$, by Nakai criterion for real divisors,
there must exist an irreducible curve $\Gamma$ on $C^{(2)}$ such
that \[((\tau(C) + 1)x - (\delta/2)).\Gamma = 0.\] Suppose $\Gamma
\equiv (n,\gamma) + \xi$. Then the result follows.

For (b). Assume $\tau(C) > a/b$. Since $D \equiv (a,b)$ is
effective, \[\tau(C)a-bg = (\tau(C),1).(a,b) \ge 0.\]

Thus, $\tau(C) \ge \frac{bg}{a} > \frac{g}{\tau(C)}$ and hence
$\tau(C) > \sqrt{g}$. From (a) an irreducible curve $\Gamma$
computing $\tau(C)$ exists and
\[D.\Gamma = an-b\gamma g = an - bn\tau(C) < 0.\] Hence, $\Gamma$ is
a component of $\supp(D)$.

Suppose $C' \equiv (n',\gamma') + \xi'$ is one of the components.
Since $C'$ is an integral curve, we have

\begin{myequation} \label{max:comp}
\tau(C)n' - \gamma'g = (\tau(C),1).((n', \gamma')+ \xi') \ge 0
\Rightarrow \tau(C) \ge \frac{\gamma'g}{n'}.
\end{myequation}

Since $(x.C') = n'$ and $(\delta.C') = 2n' + 2\gamma'g$, then
\[ \frac{\gamma'g}{n'} = \frac{(\delta.C') - 2(x.C')}{2(x.C')}.\]

Hence, $\tau(C) = \max\limits_{C' \subset \supp (D)} \{R(C')\}$.

For the particular statement, since $\Gamma(g^1_e) \equiv (e-1,1)$,
as above $\tau(C) > \sqrt{g}$ and it must be computed by one of the
irreducible components of $\Gamma(g^1_e)$ and
\[ \tau(C) = \max\limits_{C' \subset \Gamma(g^1_e)} \{ R(C') \}.\]
\end{proof}

\begin{corollary} \label{weak:in:support}
Suppose $C$ is a smooth integral curve of genus $g \ge 2$. If there
exists an effective divisor $D$ on $C^{(2)}$ such that $D \equiv
(a,b)$ with $a,b$ are positive integers and $D^2 \le 0$, then either
$\tau(C) = \sqrt{g} = \frac{a}{b},$ or an irreducible component of
$\supp(D)$ computes $\tau(C)$.
\end{corollary}
\begin{proof}
It is clear that $D^2 \le 0$ implies $a \le b\sqrt{g}$.

For the case $\tau(C) > \sqrt{g} \ge a/b$, \eqref{in:support} (b)
implies an irreducible component of $\supp(D)$ computing $\tau(C)$.

For $\tau(C) = \sqrt{g}$, since $a \tau(C)\ge bg$ (by the
effectiveness of $D$), $a = b\sqrt{g}$ and $g$ is a perfect square.
We get
\begin{myequation} \label{root:g}
\tau(C) = \sqrt{g} = {a \over b}.
\end{myequation}

\end{proof}

\begin{remark} \label{sqrt:g}
The result of \eqref{weak:in:support} implies \eqref{kouvidakis}. To
this end, since $\Gamma(g^1_d) \equiv (d-1,1)$, then
\eqref{weak:in:support} implies either
\[\tau(C) = \sqrt{g}, d-1 = \sqrt{g}{\mathrm{~and~}} g {\mbox{~is~a~perfect~square}},\] or one of the components
of $\Gamma(g^1_d)$ computes $\tau(C)$. Since $\Gamma(g^1_d)$ is
irreducible, it itself computes \[\tau(C) = g/(d-1).\]
\end{remark}

\section{Covering curves} \label{cover:curve}
In this section, we give the sufficient conditions that a curve on
the $n$-th symmetric product of $C$ is irreducible and smooth, and
takes $C$ as a covering of degree $n$. Let us recall some
definitions and equations from \cite{ACGH} we will use in the
sequel.

Let $C$ be a smooth irreducible curve of genus $g$, $J(C)$ its
Jacobian variety and $C^{(d)}$ its $d$-fold symmetric product. From
the classical theory, fixing a point $p \in C$, we can define the
morphisms
\begin{myequation} \label{abel:maps}
\begin{split}
 u_d: C^{(d)} \ra J(C) \mbox{~where~} u_d(D)=
\lb(D-d \cdot p) \mbox{~and~} D \in C^{(d)},\\
i_{d-1}: C^{(d-1)} \ra C^{(d)} \mbox{~where~} i_{d-1}(D) = D+p
\mbox{~and~} D \in C^{(d-1)}.
\end{split}
\end{myequation}

We denote the theta divisor on $J(C)$ by $\Theta$ and its numerical
class in the N{\'{e}}ron-Severi space by $\theta$. On $C^{(d)}$ we
denote the class of $u_d^*(\theta)$ by $\theta_d$ or simply $\theta$
if no confusion in context and the class of $i_{d-1}(C^{(d-1)})$ in
$C^{(d)}$ by $x_d$ or simply $x$.

Suppose $\sigma : C^{(n-1)} \times C \ra C^{(n)}$ is the natural
surjective morphism.

Set $D(n,p) = \sigma(C^{(n-1)},p)$. Clearly,

\begin{myequation} \label{cn:divisor}
D(n,p) = i_{n-1}(C^{(n-1)}) \equiv x  \mbox{~for~any~} p \in C.
\end{myequation}

Consider the morphism
\[ \phi : C^{(d-2)} \times C \ra C^{(d)}\] defined by $\phi(D,p) = D+2p$.
The image of $\phi$ is the {\it diagonal} $\delta_d$ in $C^{(d)}$.
The \cite[p.358 Proposition 5.1]{ACGH} tells that
\begin{myequation}\label{diagonal}
\delta_d \equiv 2((d+g-1)x-\theta).
\end{myequation}

Given an algebraic cycle $Z$ in $C^{(d)}$ we define two associated
cycles
\begin{myequation} \label{cycle:oper}
\begin{split}
A_k(Z) &= \{E \in C^{(d+k)} : E-D \ge 0 \fs D \in Z\},\\
B_k(Z) &= \{E \in C^{(d-k)} : D-E \ge 0 \fs D \in Z \}.\\
\end{split}
\end{myequation}

With all the above, we compute the intersection numbers for the
curves from coverings.

Suppose $\pi:C\ra H$ is an $n$-sheeted covering of $C$ over a smooth
irreducible curve $H$ of genus $h$. Then we have a curve $\Sigma :=
\{\pi^{-1}(q): q \in H\} \subset C^{(n)}$. We have the following
intersection numbers.
\begin{lemma} \cite[p.370 D-10]{ACGH}  \label{dot:basis1}
\begin{equation*}
\begin{split}
(x.\Sigma) &= 1,\\
(\delta.\Sigma) &= 2(g-1)-2n(h-1)\\
                & = \mbox{~the~degree~of~the ramification~divisor~of~}\pi,\\
(\theta.\Sigma) &= nh.\\
\end{split}
\end{equation*}
\end{lemma}

Conversely, under some hypotheses, $C$ can be an $n$-sheeted
covering. The following lemma is from \cite[1.5]{Ci}. We modified
the conditions appropriately for our purpose.

\begin{lemma}  \label{decomposition}

For any integral curve $\gamma \subset C^{(n)}$, if it satisfies
\begin{enumerate}[(i)]

\item $\gamma \not\subset \supp(D(n,p))$ for all $p \in C$ and
\item $(x.\gamma) = 1$,
\end{enumerate}
then $\gamma$ is smooth and there is a degree $n$ covering $f :C \ra
\gamma$. In particular, $\Sigma := \{ f^{-1}(q) : q \in \gamma \}
\subset C^{n}$ is smooth and isomorphic to $\gamma$.
\end{lemma}
\begin{proof}
Fix $p \in C$. Suppose $\eta \in \gamma \cap D(n,p)$. Since $\gamma
\not\subset \supp D(n,p)$,

\[1 = (x.\gamma) = D(n,p).\gamma \ge \mult_{\eta}D(p) \cdot \mult_{\eta}
\gamma = \mult_{\eta} \gamma.\]

This is true any $p \in C$. Therefore $\gamma$ intersects $D(n,p)$
only at $\eta$ and $\gamma$ is smooth at $\eta$.

Set-theoretically, a map $f: C \ra \gamma$ can be defined by
\[f(p) = D(n,p) \cap \gamma.\]

Consider the natural morphism $\sigma : C^{(n-1)} \times C \ra
C^{(n)}$ and let $p_2 : C^{(n-1)} \times C \ra C$ denote the second
factor projection.

Since $D(n,p) = \sigma(C^{(n-1)} \times p)$ and $\sigma$ is a finite
morphism, then the inverse image $\sigma^{-1}(\gamma)$ of $\gamma$
is also a curve. By the projection formula of intersection of
cycles,
\[\sigma^{-1}(\gamma).(C^{(n-1)} \times p) = \gamma. D(p) = 1.\]

Therefore there is only one integral component of
$\sigma^{-1}(\gamma)$ dominant $C$ through $p_2$, and this morphism
is finite and of degree one. Since $C$ is smooth, that integral
component is isomorphic to $C$.

On the other side, the hypothesis $\gamma \not\subset \supp(D(n,p))$
for all $p$ guarantees that no component of $\sigma^{-1}(\gamma)$
can be contained in any fibre of $p_2$. Thus, $\sigma^{-1}(\gamma)$
is irreducible and isomorphic to $C$ through $p_2$.

Then, we can define $f$ as,

\[ f := \sigma \circ (p_2|_{\sigma^{-1}(\gamma)})^{-1}.\]

Since $\sigma$ is of degree $n$, $f$ is a degree $n$ covering.

By the definition of $\Sigma$, we have a natural birational morphism
from $\gamma$ into $\Sigma$. Since $\Sigma$ comes from a degree $n$
covering, $\Sigma$ cannot be contained in the support of $D(n,p)$
for any $p \in C$. By \eqref{dot:basis1} and the above arguments,
$\Sigma$ is smooth and hence isomorphic to $\gamma$.
\end{proof}

In particular, on $C^{(2)}$, those integral curves, not a curve $C +
a$ for any fixed $a$, of intersection number one with the fibre
class $x$ have $C$ as a double cover.

\begin{remark}
With the same result of \eqref{decomposition}, the hypotheses in
\cite[1.5]{Ci} are:
\begin{enumerate}[(a)]
\item $\gamma  \not\subset \delta_n$.
\item For any $p \in C$, $D(n,p)$ cuts $\gamma$ in one single
point.
\end{enumerate}
Actually, our conditions implies these two above. To this end,
assume to the contrary of (a) that $\gamma \subset \delta_n$. A
general point on $\gamma$ is $D + 2p_0$ where $p_0 \in C$ and $D$ is
an effective divisor on $C$ of degree $n-2$. Then, there exists an
integral curve $\alpha$ in $C^{(n-1)} \times C$ with general point
on it is $(D+p_0,p_0)$ and $\sigma(\alpha) = \gamma$. The morphism
$\sigma$ is ramified along $\alpha$. Since $\gamma \not\subset
\supp(D(n,p))$ for all $p \in C$, then $\alpha$ is dominant $C$
through $p_2$. Thus, we have
\[(\gamma.x)= \gamma. D(n,p) = \sigma^{-1}(\gamma).(C^{(n-1)} \times p) \ge 2\alpha. (C^{(n-1)} \times p)
\ge 2.\] Contradiction.

For (b), it can be proved similarly as in the proof of
\eqref{decomposition}.
\end{remark}

Again suppose $\pi:C\ra H$ is an $n$-sheeted covering of $C$ over a
smooth irreducible curve $H$ of genus $h$ and we define $\Sigma :=
\{\pi^{-1}(q): q \in H\} \subset C^{(n)}$. We would like to apply
the two operations \eqref{cycle:oper} to the curve $\Sigma$. By
\eqref{decomposition}, $H \cong \Sigma$. Thus, $H$ is embedded into
$C^{(n)}$. From the definitions, $B_{n-d}(H) = B_{n-d}(\Sigma)$ is a
$1$-cycle in $C^{(d)}$ for any $2\le d \le n$. Moreover, with
\eqref{dot:basis1} and \cite[p.368 D-2 and D-8]{ACGH}, ones have the
following intersection numbers.

\begin{lemma}
\begin{equation*}\label{dot:basis2}
\begin{split}
(x.B_{n-d}(H)) &= \left(\begin{array}{c} n-1\\n-d \end{array}\right),\\
(\theta.B_{n-d}(H)) &= nh\left(\begin{array}{c}n-2\\n-d \end{array}\right)+g\left(\begin{array}{c}n-2\\n-d-1 \end{array}\right),\\
(\delta.B_{n-d}(H)) &= 2(d+g-1)\left(\begin{array}{c} n-1\\n-d \end{array}\right)-2(nh\left(\begin{array}{c}n-2\\n-d \end{array}\right)+g\left(\begin{array}{c}n-2\\n-d-1 \end{array}\right)).\\
\end{split}
\end{equation*}
In particular, for $d=2$,
\begin{myequation}\label{dot:basis3}
\begin{split}
(x.B_{n-2}(H)) &= n-1,\\
(\delta.B_{n-2}(H)) &= 2(g-1)-2n(h-1)\\
                    &= \mbox{~the~degree~of~the ramification~divisor~of~}\pi.\\
\end{split}
\end{myequation}
\end{lemma}

\begin{corollary} \label{one:decomp}
Under the above settings, if $\Gamma$ is an irreducible component of
$B_{n-d}(H)$ for some $2 \le d \le n$ and $(x.\Gamma) = 1$, then $C$
is a $d$-sheeted cover of $\Gamma$ and $\pi$ factorizes through this
cover.
\end{corollary}
\begin{proof}
Since $C \ra H$ is a covering, $\Gamma \not\subset \supp(D(d,p))$
for all $p \in C$. From \eqref{decomposition}, $C \ra \Gamma$ is a
covering of degree $d$. By the definition of $B_{n-d}(H)$, $\pi: C
\ra H$ must factorize through $\Gamma$.
\end{proof}

\section{Characterization of double cover curves}
The Lemma \eqref{decomposition} tells the geometric meaning of
intersection number one with the fibre class $x$. This leads us a
way to characterize those double covering curves. To give the proof
of \eqref{main:result}, we need two more well-known results.

\begin{theorem} \cite[Theorem 3.5]{Ac}  \label{accola}
Let $C$, $C_1$ and $C_2$ be smooth projective integral curves of
genera $g$, $g_1$ and $g_2$ respectively. Let $\pi_1 : C  \ra C_1$
and $\pi_2 : C  \ra C_2$ be $n_1$ and $n_2$-sheeted coverings.
Assume there does not exist smooth integral curve $\Gamma$ of genus
$h < g$ with coverings $f : C \ra \Gamma$ and $\alpha_i : \Gamma \ra
C_i$ such that
\[ \pi_i = \alpha_i \circ f \mymbox{for} i=1,2.\]  Then,

\[g \le n_1g_1 + n_2g_2 + (n_1 -1)(n_2-1).\]
\end{theorem}

Let $W_d^r(C)$ denote the subvariety of $\mathrm{Pic}^d(C)$
parameterizing complete linear series of degree $d$ and dimension at
least $r$:
\[ W_d^r(C) = \{ |D| : \deg D = d \mbox{~and~} r(C) \ge r \}.\]

\begin{lemma} \cite[p.181 3.3]{ACGH}  \label{brill:noether}
Suppose $r \ge d-g$. Then every component of $W_d^r(C)$ has
dimension greater or equal to the Brill-Noether number

\[\rho(g,r,d) = g - (r+1)(g-d+r).\]
\end{lemma}

The following gives the core part of \eqref{main:result}.

\begin{proposition} \label{double:cover}
Fix an integer $k \ge 0$. Suppose $C$ is a smooth integral curve of
genus $g
> 2k+1$. Assume $C$ admits a base-point-free pencil $g^1_d$ of
degree $d \ge 2$ with $d <g-2k+3$. Furthermore, assume
$\Gamma(g^1_d)$ contains an irreducible component computing
$\tau(C)$. Then,
\begin{enumerate}

\item[]  \begin{equation*}
\begin{split}{}
 &\tau(C)  \ge g-k \mbox{~if~and~only~if~}\\
 &C \mbox{~is~a~double~
cover~of~a~smooth~integral~curve~} H \\
&\mbox{of~genus~} q
\mbox{~with~} q \le \frac{k}{2}.\\
\end{split}
 \end{equation*}

\item[] Furthermore in this case $H$ is isomorphic to the component of $\Gamma(g^1_d)$
computing $\tau(C)$ and
\[\tau(C) = g-2q.\]
\end{enumerate}

\end{proposition}

\begin{proof}
Suppose $\Sigma$ is the irreducible component of $\Gamma(g^1_d)$
computing $\tau(C)$. Then, \[\tau(C) = \frac{(\delta.\Sigma) -
2(x.\Sigma)}{2(x.\Sigma)}.\]

For any irreducible component $C'$ of $\Gamma(g^1_d)$, let us set
\[R(C') := \frac{(\delta.C') -
2(x.C')}{2(x.C')}.\] It means that $\tau(C) = R(\Sigma)$.

It is clear that $(x.\gamma) \ge 1$ for any integral curve $\gamma$
on $C^{(2)}$.

It splits into two cases.
\begin{enumerate} [(i)]
\item For $(x.C') = 1$. Since $C' \subset \Gamma(g^1_d)$ and $g^1_d$ is base-point-free,
$C' \not\subset D(2,p)$ for any fixed $p \in C$. By
\eqref{decomposition}, this induces a double covering of $C$ onto
$C'$. By \eqref{dot:basis1},
$$(\delta.C') = 2g-4g'+2$$ where $g'$ is the genus of $C'$. This
computes
$$R(C') = g-2g'.$$
\item $(x.C') \ge 2$. Since $\Gamma(g^1_d)$ comes from covering of $C$,
there is no component of $\Gamma(g^1_d)$ can be the diagonal of
$C^{(2)}$. Thus, $2d-2+2g =  (\delta.\Gamma(g^1_d)) \ge (\delta.C')
\ge 0$. Since $d< g-2k+3$, this gives
$$R(C') \le \frac{(d-1)+g -2}{2} < g-k.$$
\end{enumerate}

\begin{enumerate}
\item["$\Leftarrow$":] Suppose $C \ra H$ is a double cover of a
curve $H$ of genus $q$ for $0 \le q \le \frac{k}{2}$. Then, $H$ can
be embedded into $C^{(2)}$ and $((\tau(C)+1)x - (\delta/2)).H \ge
0$. Thus, from \eqref{dot:basis1},
\[\tau(C) \ge g-2q \ge g-k.\]
The irreducible component $\Sigma$ must satisfy (i) and hence
\[\tau(C) = g-2g'\ge g-2q \mbox{~and~}0 \le g' \le q \le \frac{k}{2}.\]
If $H$ and $C'$ are not isomorphic, then $C \ra C'$ and $C \ra H$
are two distinct double covers, by \eqref{accola},
\[g \le 2q + 2g' + 1.\] However, $2q+2g'+1 \le 2k+1 <
g$. Contradiction. Hence, $H \cong C'$ and \[\tau(C) = R(C') = g-2g'
= g-2q.\]
\item["$\Rightarrow$":] Assume $\tau(C) \ge g - k$. Again the component $\Sigma$ must be in the case (i). Then, \[\tau(C) = R(C') = g-2g'.\] Thus, $0 \le g' \le \frac{k}{2}$.
\end{enumerate}
\end{proof}

If $g$ is large enough, there always exists a base-point-free pencil
satisfying the hypotheses of \eqref{double:cover}. This is a direct
application of the Brill-Noether number.

\begin{proof}[Proof of \eqref{main:result}]
Let $d$ be the gonality of $C$, i.e. $d$ is the smallest integer
such that $C$ has a base-point-free degree $d$ pencil. By
\eqref{brill:noether},

\[2 \le d \le \left[\frac{g+3}{2}\right].\]

Let $g^1_d$ denote a base-point-free pencil of degree $d$. Since $g
> 4k-3$ and $g
> 2k+1$, we have the following inequalities respectively,
\[ d \le \frac{g+3}{2} < g-2k+3\] and

\[d-1 \le \frac{g+3}{2}-1 < g-k.\]

To apply \eqref{double:cover}, we have to check $\Gamma(g^1_d)$
contains a component computing $\tau(C)$.
\begin{enumerate}
\item["$\Rightarrow$":]

Assume $\tau(C) \ge g-k$, then $$\tau(C) > d-1.$$ Thus
\eqref{in:support} implies that a component of $\Gamma(g^1_d)$
computes $\tau(C)$. Then, \eqref{double:cover} implies the
conclusion.

\item["$\Leftarrow$":] If $C \ra H$ is a double cover of a curve $H$ of
genus $q$ for $q \le \frac{k}{2}$. Then, $H$ can be embedded into
$C^{(2)}$ and $((\tau(C)+1)x - (\delta/2)).H \ge 0$. Thus,
$$\tau(C) \ge g-2q \ge g-k.$$ Again, \eqref{in:support} implies
a component of $\Gamma(g^1_d)$ to compute $\tau(C)$. Then the
conclusion follows from \eqref{double:cover}.
\end{enumerate}
\end{proof}

\begin{remark}
\begin{enumerate}[(i)]
\item The bound $\max \{4k-3, 2k+1\} = 4k-3$ when $k \ge 2$. We
put it in the theorem for the completeness to include the
hyperelliptic cases of $2 \le g \le 5$.

\item For the bielliptic case $k = 2$, the Theorem \eqref{main:result} (b) needs $g >
5$ to conclude $C$ is a bielliptic curve from $\tau(C) = g-2$.
However, in \eqref{debarre:result2}, Debarre proved that if $g =5$
and $\tau(C) = g-2 = 3$, then $C$ is bielliptic. This is a boundary
situation we cannot use the inequality \[ d < g -2k+3\] to eliminate
the case (ii) in the proof of \eqref{double:cover} when $g=5$, $k=2$
and $d=4$. Instead, Debarre used an analysis on the orbits of the
Galois group $G$ of the base-point-free $g^1_4$ on $C$. By this way
the irreducible component of $\Gamma(g^1_4)$ computing $\tau(C)$
must have intersection number one with the fibre class $x$, i.e.
must satisfy the case (i) in the proof of \eqref{double:cover}. Then
going through the arguments follow in the proof, one has the
conclusion that $C$ is bielliptic when $g=5$ and $\tau(C) = 3$.
\end{enumerate}
\end{remark}

\section{Examples}
If an integral curve $C$ is  a cover of higher degree, the value of
$\tau(C)$ can jump down to equal or less than half its genus. We are
computing some examples to illustrate this situation.

Suppose $C$ is a smooth integral curve of genus $g$. For any
integral curve $\gamma$ on $C^{(2)}$, keep using the notation
\[R(\gamma) := \frac{(\delta.\gamma) -
2(x.\gamma)}{2(x.\gamma)}.\]

\begin{enumerate}[(I)]
%\item[(a)]
%For $g = 0,1$ or $C$ is hyperelliptic, \[\tau(C) = g.\] For $g = 3$
%and $C$ is non-hyperelliptic, then $C$ is a plane quartic and
%\[\tau(C) = \frac{9}{5}\] (\cite{De} Proposition 8)(\cite{Kon}
%Proposition 3.4)(\cite{Ro} 2.2).
\item Suppose $C$ is a trigonal of genus $g \ge4$ (i.e. non-hyperelliptic and admits a base-point-free pencil $g^1_3$).
From \eqref{dot:basis3}, \[(x.\Gamma(g^1_3)) = 2.\] If
$\Gamma(g^1_3)$ is reducible, it has at most two components. Let
$C'$ be one of them, then $(x.C') = 1$. However, \eqref{one:decomp}
implies the morphism by $g^1_3$ factorizing through $C'$. This is
impossible as $3$ is prime. Hence, $\Gamma(g^1_3)$ is irreducible.
For $g \ge 4$, by \eqref{kouvidakis},
\[\tau(C) = \frac{g}{2}.\]

\item \cite[Remark 9]{De} Suppose $C$ is a tetragonal of genus $g \ge 9$.
If $\tau(C) = \sqrt{g}$, then \eqref{root:g} implies
\[\frac{g}{d-1} = \tau(C).\] Then, $g = 9$  and $\tau(C) = 3$.

Assume $\tau(C) > \sqrt{g}$.

If $\Gamma(g^1_4)$ is irreducible, then by \eqref{kouvidakis}
\[\tau(C) = \frac{g}{3}.\]

If it is reducible, then \eqref{in:support} implies one of the
components computing $\tau(C)$. Since \[(x.\Gamma(g^1_4)) = 1 + 2 =
1 + 1+ 1,\] $\Gamma(g^1_4)$ can contain two or three components and
in each case the morphism by $g^1_4$ factorizes through a
$2$-sheeted covering onto a curve $H$ by \eqref{one:decomp}. We can
use \eqref{dot:basis1} and \eqref{in:support} to compute $\tau(C)$.

For the case "$1+2$", $H$ computes $g-2h$ and $\Gamma(g^1_4)-H$
computes $h$, thus \[\tau(C) = \max\{h,g-2h\}\] where $h$ is the
genus of $H$.

For the case "$1+1+1$", \[\tau(C) = \max\{g-2h_1,g-2h_2,g-2h_3\}\]
where $h_i$ is the genus of $H_i$ and $H_i$'s are the components of
$\Gamma(g^1_4)$ for $i=1,2,3$.

\item Suppose $C$ is a $5$-gonal of genus $g \ge 16$. Since
$5$ is prime, $\Gamma(g^1_5)$ cannot decompose to a component $C'$
with $(x.C') =1$. Since $(x.\Gamma(g^1_5)) = 4$, $\Gamma(g^1_5)$ can
only decompose into at most two components and they have
intersection number with $x$ is $2$. Let $C'$ be one of them, then
$$\frac{(\delta.C') - 2(x.C')}{2(x.C')} \le \frac{(\delta.\Gamma(g^1_5)) - 2(x.C')}{2(x.C')} \le \frac{g+2}{2}.$$ Thus,
\[\frac{g}{4} \le \tau(C) \le \frac{g+2}{2}.\]

\item Suppose $C$ is a triple cover of a curve $H$ of genus $h$
and $H$ admits a base-point-free pencil $h^1_2$ (i.e. elliptic or
hyperelliptic). The composition of the triple cover and the $h^1_2$
on $H$ is a base-point-free pencil $g^1_6$ on $C$. We would like to
find a lower bound $l$ of $g$ such that \[\tau(C) = \frac{g-3h}{2}\]
when $g \ge l$. First $g >  (6-1)^2 = 25$ such that $\tau(C) >
\sqrt{g}$. Then, by \eqref{in:support}  an irreducible component of
$\Gamma(g^1_6)$ computes $\tau(C)$. Since $g^1_6$ is compounded by
$H$, \[\Gamma(g^1_6) \supset B_1(H)\] (the definition of $B_k(Z)$ in
above). From \eqref{dot:basis3}, $(x.B_1(H)) = 2$. Using the similar
argument as for trigonals, we can see that $B_1(H)$ is irreducible.
Then, \[R(B_1(H)) = \frac{g-3h}{2}.\]

For the other components in $\Gamma(g^1_6)$, there are three cases.\\
$(x.(\Gamma(g^1_6)- B_1(H)) = 3 = 1+2 =
1+1+1$.\\
For the case "$3$" and $\Gamma(g^1_6)- B_1(H) = \Gamma$,
\[R(\Gamma) = h.\]\\
For the case "$1+2$" and $\Gamma(g^1_6)- B_1(H) = \Gamma_1 +
\Gamma_2$, we have \[R(\Gamma_1) = g-2g_1 \mbox{~and~} R(\Gamma_2)
= \frac{3h-g+2g_1}{2}\] where $C \ra \Gamma_1$ is a double cover and $g_1$ is the genus of $\Gamma_1$.\\
For the case "1+1+1" and $\Gamma(g^1_6)- B_1(H) = \Gamma_1 +
\Gamma_2 + \Gamma_3$, we have \[R(\Gamma_i) = g-2g_i\] where $C \ra
\Gamma_i$ is a double cover and $g_i$ is the genus of $\Gamma_i$ for
$i=1,2,3$.

By \eqref{accola}, we have \[g-2g_i \le 2+3h\] for those double
cover components in all cases.

By Hurwitz's formula, \[g+1 \ge 2g_i.\] This implies in the case
"$1+2$",

\[R(\Gamma_2) \le  \frac{3h+1}{2}.\]

Thus, to make \[R(B_1(H)) \ge \max\{h, 2+3h, \frac{3h+1}{2}\}\] it
requires $g \ge 9h+4$. Hence, $l= \max\{9h+4,26\}$.

In particular, if $C$ is a triple cover of an elliptic curve, then

$$ \tau(C) = \frac{g-3}{2} \mbox{~if~} g \ge 26.$$
\end{enumerate}

The critical point to prove \eqref{double:cover} is
\eqref{decomposition}, which tells us very well that when
$(x.\Gamma) = 1$ we have a double cover from $C$ onto $\Gamma$.
However, we cannot control those components with $(x.\Gamma) \ge 2$.
Although we know they might come from a higher degree covering of
$C$ onto some curve \eqref{dot:basis2}, the intersection number is
an necessary condition but not sufficient. For instance, the example
(III) above, it is possible that $B_3(\Gamma(g^1_5)$ has a component
$\Gamma$ with $(x.\Gamma) = 2$ but $g^1_5$ cannot be factorized.

%
% BibTeX users please use
\bibliographystyle{alpha}
\bibliography{references}
\end{document}